\RequirePackage[l2tabu, orthodox]{nag}

\documentclass[12pt]{amsart}
\usepackage{fullpage,url,amssymb,enumerate,colonequals}
\usepackage[all]{xy} 
\usepackage{mathrsfs} 

\usepackage[OT2,T1]{fontenc}
\DeclareSymbolFont{cyrletters}{OT2}{wncyr}{m}{n}
\DeclareMathSymbol{\Sha}{\mathalpha}{cyrletters}{"58}

\usepackage{color}

\newcommand{\defi}[1]{\textsf{#1}} 

\newcommand{\Aff}{\mathbb{A}}
\newcommand{\C}{\mathbb{C}}
\newcommand{\F}{\mathbb{F}}

\newcommand{\PP}{\mathbb{P}}
\newcommand{\Q}{\mathbb{Q}}
\newcommand{\R}{\mathbb{R}}
\newcommand{\Z}{\mathbb{Z}}
\newcommand{\Qbar}{{\overline{\Q}}}

\newcommand{\Zbar}{{\overline{\Z}}}
\newcommand{\kbar}{{\overline{k}}}

\newcommand{\Fbar}{{\overline{\F}}}

\newcommand{\mm}{\mathfrak{m}}

\newcommand{\OO}{\mathscr{O}}
\newcommand{\VV}{\mathscr{V}}

\DeclareMathOperator{\Aut}{Aut}

\DeclareMathOperator{\Frob}{Frob}
\DeclareMathOperator{\Gal}{Gal}

\DeclareMathOperator{\KP}{KP}

\DeclareMathOperator{\rank}{rank}

\newcommand{\et}{{\operatorname{et}}}

\newcommand{\GL}{\operatorname{GL}}
\newcommand{\HH}{{\operatorname{H}}}

\newcommand{\M}{\operatorname{M}}

\newcommand{\del}{\partial}
\newcommand{\directsum}{\oplus} 
\newcommand{\injects}{\hookrightarrow}
\newcommand{\isom}{\simeq}

\newcommand{\surjects}{\twoheadrightarrow}
\newcommand{\tensor}{\otimes} 
\newcommand{\To}{\longrightarrow}
\newcommand{\union}{\cup} 

\newtheorem{theorem}{Theorem}[section]

\newtheorem*{mordellconjecture}{Mordell conjecture}

\theoremstyle{definition}

\theoremstyle{remark}
\newtheorem{remark}[theorem]{Remark}

\newtheorem{example}[theorem]{Example}

\makeatletter
\g@addto@macro\bfseries{\boldmath} 
\makeatother

\usepackage{microtype}
\usepackage{tikz}

\usepackage[
	pagebackref,
	pdfauthor={Bjorn Poonen}, 
	pdftitle={A p-adic approach to rational points on curves},
]{hyperref}
\usepackage[alphabetic,backrefs,lite,nobysame]{amsrefs} 

\begin{document}

\title{A \MakeLowercase{$p$}-adic approach to rational points on curves}
\author{Bjorn Poonen}
\thanks{This article is associated with a lecture given January 17, 2020 in the Current Events Bulletin at the 2020 Joint Mathematics Meetings in Denver.  The writing of this article was supported in part by National Science Foundation grant DMS-1601946 and Simons Foundation grants \#402472 (to Bjorn Poonen) and \#550033.}
\address{Department of Mathematics, Massachusetts Institute of Technology, Cambridge, MA 02139-4307, USA}
\email{poonen@math.mit.edu}
\urladdr{\url{http://math.mit.edu/~poonen/}}
\date{June 1, 2020}

\begin{abstract}
In 1922, Mordell conjectured the striking statement 
that for a polynomial equation $f(x,y)=0$,
if the topology of the set of complex number solutions is complicated enough,
then the set of rational number solutions is finite.
This was proved by Faltings in 1983, and again by a different method 
by Vojta in 1991, but neither proof provided a way to provably find all
the rational solutions, so the search for other proofs has continued.
Recently, Lawrence and Venkatesh found a third proof, 
relying on variation in families of $p$-adic Galois representations;
this is the subject of the present exposition.
\end{abstract}

\maketitle

\section{The Mordell conjecture}\label{S:Mordell}

\subsection{Rational points on curves}

The equation $x^2+y^2=z^2$ has infinitely many solutions in 
integers satisfying $\gcd(x,y,z)=1$.
Equivalently, the circle $x^2+y^2=1$ has infinitely many \emph{rational} points
($(3/5,4/5)$, $(5/13,12/13)$, etc.)
This can be understood geometrically: 
each line through $(-1,0)$ with rational slope
intersects the circle at one other point, which must have rational
coordinates since finding its coordinates amounts to solving a quadratic
equation over $\Q$ for which one rational root is already known.
The same argument shows that any nonsingular conic section
defined by a polynomial with rational coefficients having 
one rational point has infinitely many.

\begin{center}
\begin{tikzpicture}[scale=1.5]


  \draw [help lines] (-1.3,0) -- (1.3,0);
  \draw [help lines] (0,-1.3) -- (0,1.3);

  \draw[thick] (0,0) circle [radius=1];

  \node[below] at (0,-1.3) {\footnotesize $x^2+y^2=1$};

  \draw[fill] (-1.000,0) circle [radius=0.04];
  \draw[fill] (0,-1.000) circle [radius=0.04];
  \draw[fill] (0,1.000) circle [radius=0.04];
  \draw[fill] (1.000,0) circle [radius=0.04];
  \draw[fill] (-0.8000,-0.6000) circle [radius=0.04];
  \draw[fill] (-0.8000,0.6000) circle [radius=0.04];
  \draw[fill] (-0.6000,-0.8000) circle [radius=0.04];
  \draw[fill] (-0.6000,0.8000) circle [radius=0.04];
  \draw[fill] (0.6000,-0.8000) circle [radius=0.04];
  \draw[fill] (0.6000,0.8000) circle [radius=0.04];
  \draw[fill] (0.8000,-0.6000) circle [radius=0.04];
  \draw[fill] (0.8000,0.6000) circle [radius=0.04];
  \draw[fill] (-0.9231,-0.3846) circle [radius=0.04];
  \draw[fill] (-0.9231,0.3846) circle [radius=0.04];
  \draw[fill] (-0.3846,-0.9231) circle [radius=0.04];
  \draw[fill] (-0.3846,0.9231) circle [radius=0.04];
  \draw[fill] (0.3846,-0.9231) circle [radius=0.04];
  \draw[fill] (0.3846,0.9231) circle [radius=0.04];
  \draw[fill] (0.9231,-0.3846) circle [radius=0.04];
  \draw[fill] (0.9231,0.3846) circle [radius=0.04];
  \draw[fill] (-0.8824,-0.4706) circle [radius=0.04];
  \draw[fill] (-0.8824,0.4706) circle [radius=0.04];
  \draw[fill] (-0.4706,-0.8824) circle [radius=0.04];
  \draw[fill] (-0.4706,0.8824) circle [radius=0.04];
  \draw[fill] (0.4706,-0.8824) circle [radius=0.04];
  \draw[fill] (0.4706,0.8824) circle [radius=0.04];
  \draw[fill] (0.8824,-0.4706) circle [radius=0.04];
  \draw[fill] (0.8824,0.4706) circle [radius=0.04];
  \draw[fill] (-0.9600,-0.2800) circle [radius=0.04];
  \draw[fill] (-0.9600,0.2800) circle [radius=0.04];
  \draw[fill] (-0.2800,-0.9600) circle [radius=0.04];
  \draw[fill] (-0.2800,0.9600) circle [radius=0.04];
  \draw[fill] (0.2800,-0.9600) circle [radius=0.04];
  \draw[fill] (0.2800,0.9600) circle [radius=0.04];
  \draw[fill] (0.9600,-0.2800) circle [radius=0.04];
  \draw[fill] (0.9600,0.2800) circle [radius=0.04];
  \draw[fill] (-0.7241,-0.6897) circle [radius=0.04];
  \draw[fill] (-0.7241,0.6897) circle [radius=0.04];
  \draw[fill] (-0.6897,-0.7241) circle [radius=0.04];
  \draw[fill] (-0.6897,0.7241) circle [radius=0.04];
  \draw[fill] (0.6897,-0.7241) circle [radius=0.04];
  \draw[fill] (0.6897,0.7241) circle [radius=0.04];
  \draw[fill] (0.7241,-0.6897) circle [radius=0.04];
  \draw[fill] (0.7241,0.6897) circle [radius=0.04];
  \draw[fill] (-0.9459,-0.3243) circle [radius=0.04];
  \draw[fill] (-0.9459,0.3243) circle [radius=0.04];
  \draw[fill] (-0.3243,-0.9459) circle [radius=0.04];
  \draw[fill] (-0.3243,0.9459) circle [radius=0.04];
  \draw[fill] (0.3243,-0.9459) circle [radius=0.04];
  \draw[fill] (0.3243,0.9459) circle [radius=0.04];
  \draw[fill] (0.9459,-0.3243) circle [radius=0.04];
  \draw[fill] (0.9459,0.3243) circle [radius=0.04];
  \draw (-1.17889,-0.357771) -- (-0.421115,1.15777);
  \draw (-1.28284,-0.282843) -- (0.282843,1.28284);
  \draw (-1.35777,-0.178885) -- (0.957771,0.978885);
  \draw (-1.37947,-0.126491) -- (1.17947,0.726491);
  \draw (-1.39223,0.0784465) -- (1.31531,-0.463062);
  \draw (-1.4,0.) -- (1.4,0.);
  \draw (-1.39223,-0.0784465) -- (1.31531,0.463062);
  \draw (-1.37947,0.126491) -- (1.17947,-0.726491);
  \draw (-1.35777,0.178885) -- (0.957771,-0.978885);
  \draw (-1.28284,0.282843) -- (0.282843,-1.28284);
  \draw (-1.17889,0.357771) -- (-0.421115,-1.15777);
  \draw[fill] (5,0) circle [radius=0.04];
  \draw[fill] (4,1) circle [radius=0.04];
  \draw[fill] (3,0) circle [radius=0.04];
  \draw[fill] (4,-1) circle [radius=0.04];


  \draw [help lines] (4-1.3,0) -- (4+1.3,0);
  \draw [help lines] (4,-1.3) -- (4,1.3);

  \draw[thick, domain=-0.9:0.9, samples=100, smooth, variable=\t]
     plot (4 + \t, {sqrt(sqrt(1-\t*\t*\t*\t)))});
  \draw[thick, domain=-0.9:0.9, samples=100, smooth, variable=\t]
     plot (4 + \t, {-sqrt(sqrt(1-\t*\t*\t*\t)))});
  \draw[thick, domain=-0.9:0.9, samples=100, smooth, variable=\t]
     plot ({4+sqrt(sqrt(1-\t*\t*\t*\t)))} , \t);
  \draw[thick, domain=-0.9:0.9, samples=100, smooth, variable=\t]
     plot ({4-sqrt(sqrt(1-\t*\t*\t*\t)))} , \t);

     \node[below] at (4 + 0,-1.3) {\footnotesize $x^4+y^4=1$};

\end{tikzpicture}
\end{center}

In contrast, Fermat proved that the equation $x^4+y^4=z^4$ 
has no positive integer solutions.
Equivalently, the set of rational points on the plane curve $x^4+y^4=1$ 
is $\{(\pm1,0),(0,\pm1)\}$.
What about $x^4+y^4=17$?
It turns out that it too has only finitely many rational points.
(They are $(\pm 2,\pm 1)$ and $(\pm 1,\pm 2)$ \cite{Flynn-Wetherell2001}.)
More generally, for any fixed $d \ge 4$ and $a \in \Q^\times$,
the curve $x^d + y^d = a$ has only finitely many rational points.
All these finiteness claims are instances of the \emph{Mordell conjecture}, 
which states that a ``complicated enough'' curve 
has only finitely many rational points, if any at all.

In the previous paragraph, 
the condition $d \ge 4$ is what made the curve ``complicated enough''.
To state the Mordell conjecture fully, however, 
we need to consider also curves defined by several polynomials
in higher-dimensional space and to introduce the notion of genus
to measure their geometric complexity.

\subsection{Projective space}

Let $k$ be a field and let $n \in \Z_{\ge 0}$.
The set of $k$-points on $n$-dimensional \defi{affine space} 
is $\Aff^n(k) \colonequals k^n$.

Define an equivalence relation $\sim$ on $k^{n+1}-\{\vec{0}\}$ 
by $(a_0,\ldots,a_n) \sim (\lambda a_0,\ldots,\lambda a_n)$ 
for all $\lambda \in k^\times$.
Let $(a_0:\ldots:a_n)$ denote the equivalence class of $(a_0,\ldots,a_n)$.
The set of all such equivalence classes is the set
\[
	\PP^n(k) \colonequals \frac{k^{n+1} - \{\vec{0}\}}{k^\times}
\]
of $k$-points on $n$-dimensional \defi{projective space}.

The points $(a_0:\ldots:a_n) \in \PP^n(k)$ with $a_0 \ne 0$
have a unique representative of the form $(1,a_1,\ldots,a_n)$,
so they form a copy of $\Aff^n(k)$.
For each $i$, the same holds for the points with $a_i \ne 0$.
Moreover, $\PP^n(k)$ 
is the union of these $n+1$ overlapping copies of $\Aff^n(k)$.

One advantage of projective space over affine space
is that $\PP^n(\R)$ is compact
for the topology coming from the Euclidean topology on each $\R^n$;
similarly, $\PP^n(\C)$ is compact.
Related to this is that intersection theory works better in projective space:
for example, two distinct lines in $\PP^2(k)$ always meet in exactly one point.

\subsection{Projective varieties}

A finite list of polynomials $f_1,\ldots,f_m \in k[x_1,\ldots,x_n]$
defines an \defi{affine variety}\footnote{Some people require a variety to satisfy additional conditions, such as not being a union of two strictly smaller such varieties.} 
$X \subset \Aff^n$ whose set of $k$-points is
\[
	X(k) \colonequals \{\vec{a} \in \Aff^n(k) : f_1(\vec{a})=\cdots=f_m(\vec{a})=0 \}.
\]

But for a point $(a_0:\ldots:a_n) \in \PP^n(k)$,
a polynomial condition $f(\vec{a})=0$ does not necessarily make sense;
to make sure that it is unchanged by scaling $\vec{a}$,
we assume that $f$ is \defi{homogeneous}, 
a sum of monomials of the same total degree, 
such as $x_0^5 x_1^2 - x_0^4 x_1^3 + 9 x_1^7$ of degree~$7$.
A finite list of homogeneous polynomials $f_1,\ldots,f_m \in k[x_0,\ldots,x_n]$
defines a \defi{projective variety} $X \subset \PP^n$ whose set of $k$-points is
\[
	X(k) \colonequals \{(a_0:\cdots:a_n) \in \PP^n(k) : f_1(\vec{a})=\cdots=f_m(\vec{a})=0 \}.
\]
The decomposition of $\PP^n$ as a union of $n+1$ copies of $\Aff^n$
restricts to express $X$ as a union of $n+1$ affine varieties 
called \defi{affine patches}.
For each $i$, dehomogenizing $f_1,\ldots,f_m$ by setting $x_i$ equal to $1$
gives polynomials cutting out the $i$th affine patch in $\Aff^n$.

\subsection{Smooth varieties}

If a variety $Y \subset \Aff^n$ is defined by $f_1,\ldots,f_{n-r}$
such that for every field extension $L \supset k$
and point $\vec{a} \in Y(L)$,
the matrix 
$\left(\left( \frac{\del f_i}{\del x_j} \right)(\vec{a})\right) \in \M_{n-r,n}(L)$
has rank $n-r$,
then call $Y$ \defi{obviously smooth of dimension $r$};
the rank condition is the same as the Jacobian criterion 
in the implicit function theorem.
More generally, any affine or projective variety $X$ is called
\defi{smooth of dimension $r$} if 
(in a sense we will not make precise) 
it can be covered by subvarieties isomorphic 
to obviously smooth varieties $Y$ as above.

If $X$ is smooth of dimension $r$ over $\R$,
then $X(\R)$ is a smooth $\R$-manifold of dimension~$r$.
The same holds if $\R$ is replaced by $\C$ in all three places.

\subsection{Genus of a curve}

{}From now on, we consider a smooth projective curve $X$ over $\Q$,
that is, a projective variety over $\Q$ that is smooth of dimension~$1$.
We assume, moreover, that $X$ is \defi{geometrically connected},
meaning that the variety defined by the same polynomials
over an algebraically closed extension field (such as $\C$)
is nonempty and not the disjoint union of two strictly smaller varieties.
Then $X(\C)$ is a compact connected $1$-dimensional $\C$-manifold,
that is, a compact Riemann surface.
Forgetting the complex structure,
we find that $X(\C)$ is a compact connected
oriented $2$-dimensional real manifold; 
by the classification of such, 
$X(\C)$ is homeomorphic to a $g$-holed torus for some $g \in \Z_{\ge 0}$.
The integer $g$ is called the \defi{genus} of $X$.
It measures the geometric complexity of $X$.

\begin{remark}
It turns out that $g$ also equals the dimension of the space
of holomorphic $1$-forms on $X(\C)$.
One can also define $g$ algebraically, 
either by using K\"ahler differentials in place of holomorphic forms, 
or by computing the dimension of a sheaf cohomology group $\HH^1(X,\OO_X)$.
\end{remark}

\begin{example}[The Riemann sphere]
If $X=\PP^1$, then the space $X(\C) = \PP^1(\C) = \C \union \{\infty\}$
is homeomorphic to a sphere via (the inverse of) stereographic projection.
Thus $g=0$.
\end{example}

\begin{example}[Plane curves]
\label{E:plane curves}
If $X \subset \PP^2$ is a smooth projective curve defined by a degree~$d$
homogeneous polynomial, then it turns out that $g = (d-1)(d-2)/2$.
\end{example}

\begin{example}[Conic sections]
A nondegenerate conic section is a smooth curve of degree~$2$ in $\PP^2$.
By Example~\ref{E:plane curves}, it is of genus $0$.
\end{example}

\begin{example}[Elliptic curves]
An \defi{elliptic curve} is a smooth degree~$3$ curve 
$y^2 z - x^3 - Axz^2 - Bz^3=0$ in $\PP^2$ for some numbers $A,B \in \Q$.
(Dehomogenizing by setting $z=1$ gives the equation
$y^2=x^3+Ax+B$ for one affine patch.)
By Example~\ref{E:plane curves}, an elliptic curve is of genus~$1$.
\end{example}

\begin{example}[Hyperelliptic curves]
Let $f(x) \in \Q[x]$ be a nonconstant polynomial with no repeated factors.
Then $y^2=f(x)$ defines a smooth curve in $\Aff^2$.
It is isomorphic to an affine patch 
of some smooth projective geometrically connected curve $X$.
If $f$ has degree $2g+1$ or $2g+2$, then the genus of $X$ is $g$.
\end{example}

\begin{remark}
The problem of determining the rational points on a general curve
can be reduced to the problem for 
a smooth projective geometrically connected curve 
(cf.~\cite{Poonen2017-Qpoints}*{Remark~2.3.27}).
That is why it suffices to consider only the latter.
\end{remark}

\subsection{The conjecture}

\begin{mordellconjecture}[\cite{Mordell1922}, first proved in \cite{Faltings1983}]
Let $X$ be a smooth projective geometrically connected curve 
of genus $g$ over $\Q$.
If $g>1$, then $X(\Q)$ is finite.
\end{mordellconjecture}

\begin{remark}
One can say qualitatively what happens for curves of genus $0$ and $1$
as well:

\bigskip

\begin{center}
\begin{tabular}{c||c|c}
Genus $g$ & $X(\Q)$ & Some examples \\ \hline
$0$ & infinite, if nonempty & lines and conics\footnotemark\\
$1$ & can be finite or infinite & elliptic curves, \lefteqn{\dots} \\
$>1$ & finite & plane curves of degree $\ge 4$, \lefteqn{\dots}
\end{tabular}
\footnotetext{In fact, every genus~$0$ curve is isomorphic to one of these.} 
\end{center}
\end{remark}

\bigskip

Several proofs of the Mordell conjecture are known, none of them easy:
\begin{itemize}
\item 
Faltings~\cite{Faltings1983} proved the conjecture in 1983 
using methods from Arakelov theory,
a kind of arithmetic intersection theory
that combines number-theoretic data with complex-analytic data.
\item 
Vojta~\cite{Vojta1991} gave a completely different proof 
based on diophantine approximation,
a theory whose original goal was to quantify how closely
irrational algebraic numbers such as $\sqrt[3]{2}$
could be approximated by rational numbers with denominator 
of at most a certain size.
For a more elementary variant of Vojta's proof due to Bombieri,
see \cite{Bombieri1990} or \cite{Hindry-Silverman2000}.
\item
Lawrence and Venkatesh \cite{Lawrence-Venkatesh-v3} 
recently gave yet another proof.
Their proof shares some ingredients with Faltings's 
but replaces the most difficult steps 
by arguments involving $p$-adic Hodge theory.
The rest of this article is devoted to explaining some of the ideas
underlying their proof.
\end{itemize}

\begin{remark}
All of these proofs generalize to the case of curves
defined over number fields instead of just $\Q$.
(A \defi{number field} is a finite field extension over $\Q$,
such as $\Q(\sqrt{5})$.)
\end{remark}

\begin{remark}
Although the Lawrence--Venkatesh proof 
is the first \emph{complete} proof of the Mordell conjecture
using $p$-adic methods,
older $p$-adic approaches have given partial results.
Chabauty~\cite{Chabauty1941} gave a proof 
for $X$ satisfying an additional hypothesis, 
namely $\rank J(\Q) < g$ for a certain 
projective group variety $J$ associated to $X$,
the \defi{Jacobian}.
More recently, Kim~\cites{Kim2005,Kim2009} 
proposed a sophisticated extension of Chabauty's ideas, 
using the nilpotent fundamental group of $X$
as a substitute for $J$.
He proved that his approach combined with well-known conjectures 
would imply the Mordell conjecture.
Kim's approach has already led to the explicit determination of $X(\Q)$ 
for some $X$ outside the reach 
of previous methods \cite{Balakrishnan-Dogra-Mueller-Tuitman-Vonk2019},
and it may be that Kim's approach succeeds for every $X$ of genus $>1$.
\end{remark}

\begin{remark}
All the proofs so far are ineffective:
they do not prove that there is an algorithm
that takes as input the list of polynomials defining a curve $X$ of genus $>1$
and outputs the list of all rational points on $X$.
At best they give a computable upper bound on $\#X(\Q)$ in terms of $X$.
See \cite{Poonen2002-millennium} for more about the algorithmic problem.
\end{remark}

\section{Overall strategy of the Lawrence--Venkatesh proof}\label{S:strategy}

Here let us outline the strategy of Lawrence and Venkatesh,
while postponing definitions and details to later sections.

Let $X$ be a smooth projective geometrically connected curve
of genus $>1$ over $\Q$.
Lawrence and Venkatesh use two maps of sets
\begin{equation}
\label{E:two maps}
	X(\Q) \stackrel{\KP}\To \{\textup{curves}\} \stackrel{\HH^1_\et}\To \{\textup{$\Q_p$-representations of $G_\Q$}\},
\end{equation}
where each of the last two sets is really a set of isomorphism classes.
\begin{itemize}
\item 
The map $\KP$ sends a rational point $x \in X(\Q)$ 
to a curve $Y_x$ over $\Q$;
the curves $Y_x$ are the fibers of a surjective morphism $Y \to X$ 
for some \emph{2-dimensional} variety $Y$
defined in Section~\ref{S:family}.
(\defi{Fiber} means the inverse image of a point.  
$\KP$ stands for Kodaira and Parshin, 
who constructed certain $Y \to X$ for studying 
the Mordell conjecture \cite{Parshin1971}.)
\item
The map $\HH^1_\et$ sends each curve to its \'etale cohomology;
see Section~\ref{S:Galois}.
\end{itemize}
Let $\VV$ be the composition of the two maps.
To complete the proof that $X(\Q)$ is finite, 
Lawrence and Venkatesh prove that $\VV$ has finite image and finite fibers; 
see Section~\ref{S:proof}.

\section{A family of curves}
\label{S:family}

In this section, we construct the algebraic family of curves $Y \to X$.

\subsection{Fundamental group of a punctured Riemann surface}

For now, let $X$ be a compact Riemann surface of genus $g$.
Because $X$ is homeomorphic to a $4g$-gon with edges glued appropriately,
the Seifert--van Kampen theorem implies that
the fundamental group of $X$ (with respect to any base point) 
has a presentation
\[
	\pi_1(X) \isom \; \Bigr\langle a_1,b_1,\ldots,a_g,b_g \; \Bigr\rvert \; [a_1,b_1]\cdots[a_g,b_g] \Bigr\rangle,
\]
where $[a,b] \colonequals aba^{-1}b^{-1}$;
that is, $\pi_1(X)$ is the quotient of a free group on $2g$ generators
by the smallest normal subgroup containing 
the indicated product of $g$ commutators.
More generally, if $B$ is a finite subset of $X$ of size $r$, then 
\[
	\pi_1(X-B) \isom \; \Bigr\langle a_1,b_1,\ldots,a_g,b_g,c_1,\ldots,c_r \; \Bigr\rvert \; [a_1,b_1]\cdots[a_g,b_g] c_1 \cdots c_r \Bigr\rangle.
\]

\subsection{Analytic construction of a family of ramified covers}

Now fix $X$ and a finite group $G$.
Let $x \in X$.
A surjective homomorphism $\pi_1(X-\{x\}) \stackrel{\alpha}\surjects G$
defines a finite covering space of $X-\{x\}$, 
and it can be completed to 
a finite \emph{ramified} covering $Y_{x,\alpha} \to X$, 
with some branches possibly coming together above $x \in X$.

This covering depends on $\alpha$, 
but there are only finitely many $\alpha$ 
since $\pi_1(X-\{x\})$ is finitely generated.
To obtain a space not depending on a choice of any one $\alpha$,
define the finite disjoint union $Y_x \colonequals \coprod_\alpha Y_{x,\alpha}$,
which is a disconnected ramified covering of $X$.\footnote{Lawrence and Venkatesh use a variant in which $G$ has trivial center and the disjoint union is over 
\emph{conjugacy classes} of surjective homomorphisms $\alpha$;
this makes sense since the isomorphism type of $Y_{x,\alpha}$
depends only on the conjugacy class.}
As $x$ varies, the $Y_x$ vary continuously in a family.
The total space of this family
is a 2-dimensional compact complex manifold $Y$
with a proper submersion $\phi \colon Y \to X$ 
such that $\phi^{-1}(x) = Y_x$ for each $x \in X$.

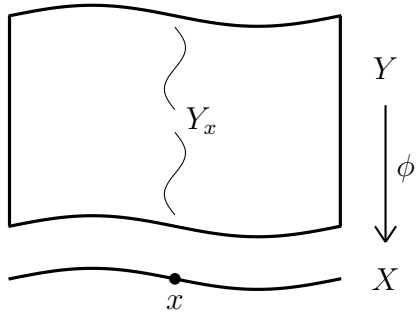
\begin{center}
\begin{tikzpicture}[scale=0.7]

  \draw[very thick, domain=-3.16:3.16, samples=100, smooth, variable=\t]
     plot (\t, {-0.2 * sin(\t r)} );

  \draw[very thick, domain=-3.16:3.16, samples=100, smooth, variable=\t]
     plot (\t, {1 -0.2 * sin(\t r)} );

  \draw[very thick, domain=-3.16:3.16, samples=100, smooth, variable=\t]
     plot (\t, {5 -0.2 * sin(\t r)} );

  \draw [very thick] (-3.14159,1) -- (-3.14159,5);
  \draw [very thick] (3.14159,1) -- (3.14159,5);

  \draw[domain=-3.14159:3.14159, samples=100, smooth, variable=\t]
     plot ({0.2 * sin(\t r)} , {2 + \t/4} );

  \draw[domain=-3.14159:3.14159, samples=100, smooth, variable=\t]
     plot ({0.2 * sin(\t r)} , {4 + \t/4} );

  \draw[fill] (0,0) circle [radius=0.1];

  \node[below] at (0,-0.1) {$x$};
  \node at (4,4) {$Y$};
  \node at (4,0) {$X$};
  \node[right] at (0,3) {$Y_x$};

  \draw [thick] (4, 3.3) -- (4, 0.7);
  \draw [thick] (3.8, 1.0) -- (4, 0.7);
  \draw [thick] (4.2, 1.0) -- (4, 0.7);
  \node[right] at (4,2.1) {$\phi$};

\end{tikzpicture}
\end{center}

\subsection{An algebraic family of curves}

The constructions above can be made algebraic, in the following sense.
Suppose that $X$ is a smooth projective connected curve over $\C$.
Then by the Riemann existence theorem,
$Y_{x,\alpha} \to X$ arises from an algebraic morphism of algebraic curves.
Moreover, there is a $2$-dimensional variety $Y$
with a morphism $\phi \colon Y \to X$ 
whose fibers are the disconnected curves $Y_x = \coprod_\alpha Y_{x,\alpha}$.

Even better, the construction is canonical enough
that if $X$ is defined over $\Q$,
then $\phi \colon Y \to X$ can be defined over $\Q$.
This is called a \defi{Kodaira--Parshin family};
see \cite{Lawrence-Venkatesh-v3}*{\S7} for details.

\begin{remark}
The curve $X$ is playing two roles: it is the base of the family $Y \to X$,
but also each fiber $Y_x$ is a ramified covering of $X$.
\end{remark}

\section{Galois representations}\label{S:Galois}

The Lawrence--Venkatesh proof makes essential use 
of $p$-adic Galois representations.
Therefore, in this section we define $\Q_p$, 
define the absolute Galois group of a field,
and give examples and properties of 
$\Q_p$-representations of the absolute Galois group of $\Q$.

\subsection{The field of \texorpdfstring{$p$}{p}-adic numbers}

Let $p$ be a prime number.
\defi{The ring of $p$-adic integers} is the inverse limit 
$\Z_p \colonequals \varprojlim \Z/p^n\Z$.
Thus an element of $\Z_p$ is a sequence $(a_1,a_2,\ldots)$
where the elements $a_n \in \Z/p^n\Z$ 
are compatible in the sense that
the natural homomorphism $\Z/p^{n+1}\Z \surjects \Z/p^n \Z$
maps $a_{n+1}$ to $a_n$ for each $n$.
For example,
\[
	(3 \bmod 5, \; 13 \bmod 5^2, \; 38 \bmod 5^3, \; \ldots) \in \Z_5.
\]
As a ring, $\Z_p$ is a domain of characteristic~$0$.
Its fraction field, denoted $\Q_p$, 
is called \defi{the field of $p$-adic numbers}.

For each $n \ge 1$, the homomorphism $\pi_n \colon \Z_p \to \Z/p^n\Z$ 
sending $(a_1,a_2,\ldots)$ to $a_n$
is surjective with kernel $p^n \Z_p$.
The kernel of $\pi_1 \colon \Z_p \to \Z/p\Z = \F_p$
is the unique maximal ideal $p\Z_p$ of $\Z_p$.
The collection of subsets $\pi_n^{-1}(a)$ 
for all $n \in \Z_{\ge 1}$ and $a \in \Z/p^n\Z$
is a basis of a topology on $\Z_p$.
Equip $\Q_p$ with the unique topology making it a topological group
having $\Z_p$ as an open subgroup.

\begin{remark}
Here we explain an alternative construction 
of $\Z_p$ and $\Q_p$ and their topologies,
producing the same results.
The \defi{$p$-adic absolute value on $\Q$}
is characterized by $\left| p^n \frac{a}{b} \right|_p \colonequals p^{-n}$
whenever $a,b,n \in \Z$ and $p \nmid a,b$;
thus a rational number is $p$-adically small 
if its numerator is divisible by a large power of $p$.
Define $\Q_p$ as the completion of $\Q$ with respect to $|\;|_p$,
just as $\R$ is the completion of $\Q$ with respect to 
the standard absolute value.
Then $|\;|_p$ extends to an absolute value on $\Q_p$.
Define $\Z_p$ as the closed unit disk $\{ x \in \Q_p : |x|_p \le 1 \}$.
Finally, $|\;|_p$ induces a metric on $\Q_p$,
which defines a topology on $\Z_p$ and $\Q_p$.
\end{remark}

Working with $\Z_p$ or $\Q_p$ amounts to working with 
infinitely many congruences at once,
but passing to the limit has advantages.
One is that one can work over a domain or field of characteristic~$0$.
Another is that some ideas from analysis over $\R$ have analogues for $\Q_p$.

Whereas number fields such as $\Q$ are examples of what are called
\defi{global fields}, $\Q_p$ is an example of a \defi{local field}.
For a more detailed introduction to $p$-adic numbers, see \cite{Koblitz1984}.

\subsection{The absolute Galois group of \texorpdfstring{$\Q$}{Q}}

A complex number is \defi{algebraic} over $\Q$
if it is a zero of some nonzero polynomial in $\Q[x]$.
The set of all algebraic numbers is a subfield $\Qbar$ of $\C$,
called an \defi{algebraic closure} of $\Q$.

Now let $K$ be a subfield of $\Qbar$.
Call $K \supset \Q$ a \defi{finite extension} if $\dim_\Q K$ is finite.
Call $K \supset \Q$ a \defi{Galois extension}  
if it is generated by the set of \emph{all} zeros of some collection
of polynomials in $\Q[x]$.\footnote{For a definition that works over an arbitrary ground field $k$ instead of $\Q$, one should require each polynomial to have distinct zeros in $\kbar$.} 
For example, $\Q(\sqrt[3]{2})$ is not a Galois extension of $\Q$,
but $\Q(\sqrt[3]{2},e^{2\pi i/3} \sqrt[3]{2},e^{4\pi i/3} \sqrt[3]{2})$ is.
The field $\Qbar$ is the union of its finite Galois subextensions $K$.

For a Galois extension $K \supset \Q$,
the \defi{Galois group} $\Gal(K/\Q)$ is the set of automorphisms of $K$
that fix $\Q$ pointwise\footnote{Fixing $\Q$ pointwise is automatic; this condition becomes relevant only over other ground fields.}.
The \defi{absolute Galois group} of $\Q$ is $G_\Q \colonequals \Gal(\Qbar/\Q)$.
Each automorphism of $\Qbar$ restricts to give an automorphism 
of each finite Galois subextension $K$,
and any compatible collection of such automorphisms 
defines an automorphism of $\Qbar$, so 
\[
	G_\Q \isom \varprojlim_{\textup{finite Galois $K \supset \Q$}} \Gal(K/\Q).
\]
Just as the inverse limit $\Z_p$ had a topology,
the inverse limit $G_\Q$ has a topology.

\begin{remark}
More generally, for any field $F$, 
one can construct a field $\overline{F}$ and topological group $G_F$.
\end{remark}

\subsection{Global \texorpdfstring{$p$}{p}-adic Galois representations}

Let $V$ be a finite-dimensional $\Q_p$-vector space.
If $\dim V = r$, then $\Aut V \isom \GL_r(\Q_p)$,
which has a topology coming from the topology of $\Q_p$.
Call a $\Q_p$-linear action of $G_\Q$ on $V$ \defi{continuous}
if the homomorphism $\rho \colon G_\Q \to \Aut V$ defined by the action
is continuous.
By a \defi{$\Q_p$-representation of $G_\Q$}
we mean a finite-dimensional $\Q_p$-vector space $V$
equipped with a continuous action of $G_\Q$.
In the next few sections, 
we give examples of such representations 
arising in number theory and arithmetic geometry.

\subsection{The cyclotomic character}

Let $m$ be a positive integer.
Define
\[
	\mu_m \colonequals \{z \in \Qbar : z^m=1 \},
\]
which under multiplication is a cyclic group of order $m$.
Thus $\mu_m$ is a free $\Z/m\Z$-module of rank~$1$.
The group $G_\Q$ acts on the group $\mu_m$.

Now fix a prime $p$, and let $m$ range through the powers of $p$.
Form the inverse limit \[
	T \colonequals \varprojlim \mu_{p^n}
\]
with respect to the homomorphisms $\mu_{p^{n+1}} \surjects \mu_{p^n}$ 
sending $\zeta$ to $\zeta^p$.
Then $T$ is a free rank~$1$ module 
under the ring $\Z_p \colonequals \varprojlim \Z/p^n\Z$,
and $G_\Q$ acts on $T$.

Next let 
\[
	V \colonequals T \tensor_{\Z_p} \Q_p.
\]
Then $V$ is a $1$-dimensional $\Q_p$-vector space,
and $G_\Q$ acts on $V$.
It follows from the definitions that the action is continuous,
so $V$ is a $1$-dimensional $\Q_p$-representation of $G_\Q$;
it is called the \defi{cyclotomic character}.

\subsection{Galois representations associated to elliptic curves}

Let $E$ be an elliptic curve over $\Q$.
It turns out that $E$ is a group variety; 
in particular, there is a map of varieties $E \times E \to E$
making $E(\Qbar)$ an abelian group.
If $P \in E(\Qbar)$, 
we may use this group law to define $3P \colonequals P+P+P$ and so on.
For each $m \ge 1$, 
it turns out that the \defi{$m$-torsion subgroup} 
\[
	E[m] \colonequals \{P \in E(\Qbar): mP=0 \}
\]
is a free $\Z/m\Z$-module of rank~$2$.
Therefore the inverse limit 
\[
	T_p E \colonequals \varprojlim E[p^n]
\]
(with respect to the homomorphisms $E[p^{n+1}] \to E[p^n]$
sending $P$ to $pP$)
is a free $\Z_p$-module of rank~$2$, called a \defi{Tate module}.
Next, 
\[
	V_p E \colonequals T_p E \tensor_{\Z_p} \Q_p
\]
is a $2$-dimensional $\Q_p$-vector space.
The continuous action of $G_\Q$ on $E(\Qbar)$
induces continuous actions on $E[p^n]$, $T_p E$, and $V_p E$.
Thus $V_p E$ is a $2$-dimensional $\Q_p$-representation of $G_\Q$.

\subsection{Galois representations associated to higher-genus curves}
\label{S:higher genus}

Let $X$ be a smooth projective geometrically connected curve
of genus $g$ over $\Q$.
If $g \ne 1$, there is no group law $X \times X \to X$,
but the Jacobian $J$ of $X$ does have a group law.
The construction of $V_p E$ generalizes to produce
a $2g$-dimensional $\Q_p$-representation $V_p J$ of $G_\Q$.

\subsection{Galois representations from \'etale cohomology}
\label{S:etale cohomology}

If $X$ is a smooth projective variety over $\Q$,
and $i \in \Z_{\ge 0}$
then the \defi{\'etale cohomology group} $\HH^i(X_{\Qbar},\Q_p)$
(which we will not attempt to define here)
is a $\Q_p$-representation of $G_\Q$.

\begin{example}
If $E$ is an elliptic curve, then it turns out that
$\HH^1(E_{\Qbar},\Q_p)$ is the dual of the representation $V_p E$.
If $X$ and $J$ are as in Section~\ref{S:higher genus},
then $\HH^1(X_{\Qbar},\Q_p)$ is the dual of $V_p J$.
\end{example}

\subsection{Semisimple representations}

Let $V$ be a $\Q_p$-representation of $G_\Q$.
Call $V$ \defi{irreducible}
if $V \ne 0$ and there is no $G_\Q$-invariant subspace $W$
with $0 \subsetneq W \subsetneq V$.
Call $V$ \defi{semisimple} 
if it is a direct sum of irreducible representations.
Maschke's theorem \cite{SerreLinearRepresentations}*{\S1.4, Theorem~2} 
states that any $\C$-representation of a finite group is semisimple, 
but this is not true for $\F_p$-representations of a finite group
of order divisible by $p$,
and $\Q_p$-representations of $G_\Q$ are more like the latter in this regard: 
they need not be semisimple.

\begin{example}
Let $\chi \colon G_\Q \to \Q_p^\times$ be the cyclotomic character.
There is a homomorphism $\log_p \colon \Q_p^\times \to \Q_p$
from the multiplicative group to the additive group; 
see \cite{Koblitz1984}*{IV.2}.
Composing these yields a nontrivial continuous homomorphism
$\lambda \colon G_\Q \to \Q_p$.
Let $V \colonequals \Q_p \directsum \Q_p$, viewed as a space of column vectors.
Let each $g \in G_\Q$ act as $\begin{pmatrix} 1 & \lambda(g) \\ 0 & 1 \end{pmatrix}$ on $V$.
The only invariant subspace of $V$ is $\Q_p \directsum 0$,
so $V$ is not semisimple.
\end{example}

\subsection{The absolute Galois group of \texorpdfstring{$\Q_p$}{Qp}}

Let $\Qbar_p$ denote an algebraic closure of $\Q_p$.
The homomorphism $\Z \injects \Z_p$
extends uniquely to $\Q \injects \Q_p$ 
and non-uniquely to $\Qbar \injects \Qbar_p$;
fix one such embedding.
Define $G_{\Q_p} \colonequals \Gal(\Qbar_p/\Q_p)$.
It turns out that $\Qbar_p$ is generated by its subfields $\Qbar$ and $\Q_p$,
so the homomorphism $G_{\Q_p} \to G_\Q$ sending each $\sigma$ to $\sigma|_{\Qbar}$
is injective.
Identify $G_{\Q_p}$ with its image, which is called 
a \defi{decomposition group} of $G_\Q$.

The absolute value $|\;|_p$ on $\Q_p$ extends in a unique way to $\Qbar_p$.
Let $\Zbar_p \colonequals \{x \in \Qbar_p : |x|_p \le 1\}$;
it is a subring.
The unique maximal ideal of $\Zbar_p$ is
$\mm \colonequals \{x \in \Qbar_p : |x|_p < 1\}$,
and the quotient $\Zbar_p/\mm$
is an algebraic closure $\Fbar_p$ of $\F_p$.
Each element of $G_{\Q_p}$ preserves $|\;|_p$ 
and hence induces an automorphism of $\Zbar_p/\mm$.
Thus we obtain a homomorphism $G_{\Q_p} \to G_{\F_p}$.
It is surjective, and its kernel $I_p \subset G_{\Q_p} \subset G_\Q$ 
is called an \defi{inertia group}.
To summarize, we have a diagram
\[
\xymatrix{
1 \ar[r] & I_p \ar[r] & G_{\Q_p} \ar[r] \ar@{^{(}->}[d] & G_{\F_p} \ar[r] & 1 \\
& & G_{\Q}.
}
\]
The \defi{Frobenius automorphism} $\Frob_p \in G_{\F_p}$
is the automorphism $x \mapsto x^p$ of $\Fbar_p$;
it generates a dense subgroup of $G_{\F_p}$ since it restricts
to a generator of each finite quotient $\Gal(\F_{p^n}/\F_p)$.
Write $\Frob_p$ also for any element of $G_{\Q_p}$
mapping to $\Frob_p \in G_{\F_p}$,
or for the corresponding element of $G_\Q$.

\subsection{Local Galois representations}

Let $p$ and $q$ be primes.
(Soon we will take $q=p$.)
A \defi{$\Q_p$-representation of $G_{\Q_q}$} is 
a finite-dimensional $\Q_p$-vector space $V$ 
equipped with a continuous action of $G_{\Q_q}$.
Call $V$ \defi{unramified} if $I_q$ acts trivially on $V$;
in that case the $G_{\Q_q}$-action can be described by one matrix,
namely the automorphism $\Frob_q|_V \in \GL(V)$ given by the action of 
any $\Frob_q \in G_{\Q_q}$.
Given $w \in \Z$, call such a $V$ \defi{pure of weight $w$} 
if the characteristic polynomial of $\Frob_q|_V$
is a polynomial in $\Z[x]$ whose complex zeros have absolute value $q^{w/2}$.

\subsection{Properties of representations coming from geometry}

Now return to a global representation $V$, 
a $\Q_p$-representation of $G_\Q$.
For each prime $q$,
restricting the $G_\Q$-action to the subgroup $G_{\Q_q}$
yields a local representation $V_q$.
For a finite set $S$ of primes, call $V$ \defi{unramified outside $S$}
if $V_q$ is unramified for every $q \notin S$.
For $w \in \Z$, call such a $V$ \defi{pure of weight $w$} 
if $V_q$ is pure of weight $w$ for every $q \notin S$.
These properties were introduced
because they hold for representations ``coming from geometry'':

\begin{theorem}
\label{T:representations from geometry} 
Each representation $\HH^i(X_{\Qbar},\Q_p)$ 
as in Section~\ref{S:etale cohomology}
is unramified outside a finite set $S$
and is pure of weight $i$ $($cf.~\cite{Deligne1974}*{Th\'eor\`eme~1.6}$)$.
\end{theorem}

\begin{remark}
\label{R:good reduction}
One can say more about $S$.
The variety $X$ can be defined by polynomials with coefficients in $\Z$.
Reducing all the coefficients of the polynomials modulo $\ell$
produces polynomials defining a variety over $\F_\ell$.
For most $\ell$, this variety is again smooth;
more precisely, this holds for all primes $\ell$ outside a finite set $S_0$.
Then in Theorem~\ref{T:representations from geometry} 
one may take $S = S_0 \union \{p\}$.
\end{remark}

\subsection{Faltings's finiteness theorem for global Galois representations}

Faltings cleverly combined a few classical facts from number theory 
(Hermite's finiteness theorem and the Chebotarev density theorem)
to prove the following finiteness statement.

\begin{theorem}[cf.~\cite{Faltings1983}*{proof of Satz~5}]
\label{T:Faltings finiteness}
Fix 
a nonnegative integer $d$, 
a prime $p$,
a finite set $S$ of primes,
and an integer $w$.
Then the set of isomorphism classes of
semisimple $d$-dimensional $\Q_p$-representations of $G_\Q$
that are unramified outside $S$ and pure of weight $w$ is finite.
\end{theorem}

\section{The Lawrence--Venkatesh proof of the Mordell conjecture}
\label{S:proof}

We now flesh out the sketch we gave in Section~\ref{S:strategy},
though we will still have to gloss over many difficult arguments.

\subsection{From rational points to representations}

Let $X$ be a smooth projective geometrically connected curve 
of genus $>1$ over $\Q$.
The goal is to prove that $X(\Q)$ is finite.

Let $G$ be a finite group.
All the claims to be made in Section~\ref{S:variation}
will be true if $G$ is chosen suitably.
(Lawrence and Venkatesh take $G$ to be 
the semidirect product $\F_q \rtimes \F_q^\times$ 
for a suitable large prime $q$.)
Let $\phi \colon Y \to X$ be the Kodaira--Parshin family of curves over $X$
defined using $G$;
let $\KP$ be the map sending each $x \in X(\Q)$ 
to the smooth projective curve $Y_x \colonequals \phi^{-1}(x)$
over $\Q$.
Let $\HH^1_{\et}$ denote the map sending 
each smooth projective curve $C$ over $\Q$
to the global Galois representation $\HH^1_{\et}(C_{\Qbar},\Q_p)$.
The composition of these, as in~\eqref{E:two maps},
is a map of sets $\VV$:
\[
\xymatrix@R=0pt{
X(\Q) \ar[r]^-{\KP} \ar@(u,u)[rr]^{\VV}  & \{\textup{curves}\} \ar[r]^-{\HH^1_\et} & \{\textup{$\Q_p$-representations of $G_\Q$}\} \\
x \ar@{|->}[r] & Y_x \ar@{|->}[r] & V_x \colonequals \HH^1_\et((Y_x)_{\Qbar},\Q_p).
}
\]
Now it turns out that
\begin{itemize}
\item The representations $V_x$ are all of the same dimension.
\item They are semisimple.\footnote{The semisimplicity is actually a difficult theorem, proved by Faltings in his paper on the Mordell conjecture.  Lawrence and Venkatesh would be ``cheating'' if they used this, so instead they give an independent argument using Hodge--Tate weights to prove that $V_x$ is semisimple for all but finitely many $x \in X(\Q)$; that is sufficient for their proof of the Mordell conjecture.}
\item They are unramified outside a set $S$ that is independent of $x$,
because one can choose a set $S_0$ as in Remark~\ref{R:good reduction} 
that works for the whole family $Y \to X$.
\item They are pure of weight~$1$.
\end{itemize}
That is, the representations $V_x$ satisfy all the conditions of Faltings's
finiteness theorem (Theorem~\ref{T:Faltings finiteness}), so
\begin{center}
   the map $\VV$ has finite image!
\end{center}

To finish the proof that $X(\Q)$ itself is finite,
one needs to show that \emph{every fiber of $\VV$ is finite},
i.e.,
that the $V_x$ vary enough that there are only finitely many $x \in X(\Q)$
mapping to any given isomorphism class of representations.

\subsection{Variation in a $p$-adic family of local Galois representations}
\label{S:variation}

The plan is to show that the global representations $V_x$ vary enough 
by showing that even their restrictions to $G_{\Q_p}$ vary enough.
These restrictions are local Galois representations indexed by $x \in X(\Q)$,
but to study them, we view them as members of a larger family
of representations, indexed by $x \in X(\Q_p)$.
Namely, for $x \in X(\Q_p)$, define the local Galois representation
\[
	V_x \colonequals \HH^1_{\et}((Y_x)_{\Qbar_p},\Q_p).
\]
Then $x \mapsto V_x$ defines the map $\VV_p$
in the following commutative diagram of sets:
\[
\xymatrix{
X(\Q) \ar@{^{(}->}[d] \ar[r]^-{\VV} & \{\textup{$\Q_p$-representations of $G_\Q$}\} \ar[d]^{\textup{restriction}} \\
X(\Q_p) \ar[r]^-{\VV_p} & \{\textup{$\Q_p$-representations of $G_{\Q_p}$}\}. \\
}
\]
To prove that $\VV$ has finite fibers,
it suffices to prove that $\VV_p$ has finite fibers.
That is, loosely speaking, one must show that the local representation $V_x$ 
varies enough as $x$ ranges over $X(\Q_p)$;
it is this claim that a large part of the Lawrence--Venkatesh article 
is devoted to.
Its proof proceeds as follows:
\begin{itemize}
\item First, $p$-adic Hodge theory relates the variation of 
the \'etale cohomology groups $V_x$ for $x \in X(\Q_p)$
to the variation of the Hodge filtration 
in the corresponding de Rham cohomology groups.
\item
The variation of the Hodge filtration 
is described by the Gauss--Manin connection,
which in down-to-earth terms 
means that it is described by the solutions to 
a system of differential equations
whose coefficients are algebraic functions on $X$ over $\Q$.
\item 
The same differential equations
describe the variation of the Hodge filtration 
for the family $Y_\C \to X_\C$ of complex projective curves.
\item
A lower bound on that variation is given by 
the monodromy of the Kodaira--Parshin family over $\C$.
\item
An extensive calculation in topology 
(involving mapping class groups, Dehn twists, and the like) 
proves that indeed the monodromy group is large enough.
\end{itemize}
This completes the proof of the Mordell conjecture.

\begin{remark}
Lawrence and Venkatesh show that their approach has applications
beyond rational points on curves.
In particular, 
using recent work of Bakker and Tsimerman~\cite{Bakker-Tsimerman2019},
they prove that certain affine varieties $F$ of higher dimension
(moduli spaces of smooth hypersurfaces in projective space)
have few \emph{integral} points, 
where ``few'' means that they are contained in a subvariety of $F$
of lower dimension.
\end{remark}



\end{document}